\documentclass[a4paper,12pt,reqno]{amsart}

\newcommand{\be}{\begin{equation}}
\newcommand{\ee}{\end{equation}}
\newcommand{\beq}{\begin{eqnarray}}
\newcommand{\eq}{\end{eqnarray}}
\newcommand{\beqq}{\begin{eqnarray*}}
\newcommand{\eqq}{\end{eqnarray*}}
                 \usepackage{graphicx}           % для включения графических изображений
\begin{document}

\begin{center}
{\bf On  $T_{2n-1}^\perp$ spaces}
\end{center}

\vskip .3cm
\begin{center}
A.G. Babenko,  Yu.V. Kryakin
\end{center}
\vskip .3cm

\centerline{\em Dedicated to Victor Kolyada on the occasion of
his  $61{}^{\text{th}}$ birthday}

\begin{center}\today\end{center}

%\vskip .3cm
%\begin{quote}
%\scriptsize \mbox{~~~}  {\it Аннотация.}
%\end{quote}

%\begin{quote}
%\scriptsize \mbox{~~~}  {\it Abstract.}

 %\end{quote}

\vskip 1cm
\begin{center}
{\bf 1 Problem}
\end{center}
\vskip .2cm

This paper is devoted to the inequalities for mean values of
functions from $T_{2n-1}^\perp$. By $T_{2n-1}$ we denote a space
of $1$--periodic real trigonometric polynomials of degree $ \le
n-1$ and  by $T_{2n-1}^\perp$ we denote its orthogonal factor in
decomposition of space\break
\mbox{${L_\infty(T)}=T_{2n-1}\bigoplus T_{2n-1}^\perp$,} \
{$T:=R/Z$}.
 The typical result is:

 $$
 g \in T_{2n-1}^\perp \implies |(g* \chi_h)(x) | \le c(h,n)\,
 {{\rm ess} \sup_{\hspace{-5mm}x\in{T}}} |g(x)|,\quad n\in N,\quad h>0,
 \eqno (1)
 $$
where  {$\chi_h$} is the characteristic function  of the
interval $(-h/2, h/2)$ normed by condition $$ \int_{-h/2}^{h/2}
\chi_h(t) \, dt =1. $$
%The inequality (1) states that the oscillation of the function from
%$T_{2n-1}^\perp$ on the interval $(x-h/2, x+h/2)$  is controlled by constant $c(h,n)$.

We describe the  procedure that  gives the sharp values of
constants $$ c(h,n):=\sup_{g \in T_{2n-1}^\perp} \frac{ \|
g*\chi_h \|}{\|g\|} $$ for all values $n$ and $h$. We present the
results of calculations for $c(h,n)$ in the nontrivial principal
case\; $h=1/n$,\; {$n\ge2$}.

The simple proof of the classical Jackson inequality  in the case
of the second modulus of continuity may be considered as the
consequence of the estimates \mbox{{$c(h,n)<1$.}} The problems of
the sharp constants in classical Stechkin's inequality are also
discussed.

\vskip .5cm
\begin{center}
{\bf 2 Notation and equivalent form of problem }
\end{center}
\vskip .2cm

Let  $C_n \ ( S_{n-1}) $ be a space of real even (odd)
trigonometric polynomials

$$
\sum_{j=0}^{n-1} a_j \cos (2\pi jx)  \quad  \left( \sum_{j=1}^{n-1} b_j
\sin (2\pi jx) \right).
$$
Denote by  $T_{2n-1}:= C_n \bigoplus S_{n-1}$ a space of all real
trigonometric polynomials.

%Let $C_n^\perp \ ( S_n^\perp, T_{2n-1}^\perp ) $ be an orthogonal
%complement of trigonometric space, in other words a space of
%1-periodic, integrable on $(-1/2,1/2)$ functions, which are orthogonal to $C_n
%\ ( S_n, T_{2n-1}) $ with respect to the scalar product
Let $C_n^\perp$ be an orthogonal complement of trigonometric space
$C_n$. In other words, a space of even 1--periodic, functions
{from $L_\infty(T)$} which are orthogonal to $C_n$ with respect to
the scalar product
 $$
 (f,g)=\int_{-1/2}^{1/2}f(t)g(t)\,dt.
 $$
One can define spaces $S_n^\perp,$ $T_{2n-1}^\perp$ in the same
manner. \noindent Central Steklov's means are the convolutions of
integrable functions with $ \chi_h (t)$: $$ S_h(f,x):= \frac{1}{h}
\int_{x-h/2}^{x+h/2} f(t) \, dt =  \frac{1}{h} \int_{-h/2}^{h/2}
f(x+t) \, dt = $$ $$ \frac{1}{h} \int_{-h/2}^{h/2} f(x-t) \, dt =
\int_{{{R}}} f(t) \chi_h (x-t) \, dt =: (f*\chi_h )(x). $$ If $$
f(x+1)=f(x), \quad \widetilde \chi_h (x):= \sum_{j=-\infty}^\infty
\chi_h(x+j), $$ then $$ (f*\chi_h )(x) = \int_{T}f(t)
\widetilde\chi_h (x-t) \, dt =:f\odot \widetilde \chi_h (x).
 $$
\vskip .2cm \noindent For $g \in T_{2n-1}^\perp, \quad \tau \in
T_{2n-1}, $

 $$
| S_h(g,x)|= |(g*\chi_h)(x)| =  | g \odot (\widetilde\chi_h-\tau)(x)| \le
\inf_\tau \|
\widetilde\chi_h - \tau \|_1 \| g \|_\infty = E_n (\widetilde \chi_h)_1 \|g \|.
 $$
%{here $\|\cdot\|_1:=\|\cdot\|_L$ and $\|\cdot\|:=\|\cdot\|_\infty:=\|\cdot\|_{L_\infty}$
%denotes norms in $L$ and $L_\infty$
%respectively.}

\noindent If  $\tau_h$ is the polynomial of the best
$L$--approximation  of the characteristic function $\widetilde\chi_h$,
then
 $$
 g_h {:=} \mbox{sign } (\widetilde\chi_h - \tau_h) \in T_{2n-1}^\perp
 \quad \mbox{(A.A. Markov (1898) \, \cite{M}),}
 $$
 and
 $$
 |(g_h*\chi_h)(0)| = |(g_h\odot \widetilde \chi_h)(0)| =  \int_{T}
|\widetilde \chi_h(t) - \tau_h (t) | \, dt = E_n (\widetilde\chi_h)_1 \|g_h \|.
 $$
\noindent So, the best constant in the inequality (1) {is}

 $$
 c(h,n) = E_n (\widetilde\chi_h)_1.
\eqno(1')
 $$

\vskip 0.5 cm
\begin{center}
{\bf 3  About the value $E_n (\widetilde\chi_h)_1$. }
\end{center}
\vskip .2cm

{Firstly note that for $h>1$
 $$
 E_n(\widetilde\chi_h)_1=\frac{\{h\}}{h}E_n(\widetilde\chi_{\{h\}})_1,  \eqno(1'')
 $$
where $\{h\}$ is the fractional part of $h$.}

It is not difficult to find $E_n (\widetilde\chi_h)_1$ for the special values {of $h$}.
% Canonical signature
{Classical signum-function}
 $$
 \mbox{sign } (\cos 2\pi nt){\in C_n^\perp}
 $$
allows us {(see \cite[Theorem 1.3.1]{bk1}, \cite[(5.1),
(5.2)]{bku})} to find the {value $E_n(\widetilde\chi_h)_1$} for
\vskip .2cm \qquad\qquad
 $
 h\in  M_n={{\left(0,\frac{1}{2n}\right]\bigcup
 \left(1-\frac{1}{2n},1\right]\bigcup_{j=2}^n\left\{\frac{2j-1}{2n}\right\}}}
 $
:

 {$$
 E_n(\widetilde\chi_h)_1 = \begin{cases}
 1, & h \in \left( 0, \frac{1}{2n} \right], \\ \\
 \frac{1}{2nh}, &  h= \frac{2j-1}{2n}, \quad j=2, \ldots,n, \\
 \\
 \frac{1-h}{h}, & h \in \left( 1-\frac{1}{2n},1 \right],
 \end{cases}
 $$}
%%%%%In the case
{and to prove that
%%%%% $$
%%%%%  E_n(\chi_h)_1 < \frac{1}{2nh},\quad
%%%%%  h\in(0,1]\setminus M_n,\
%%%%%  \eqno (5')
%%%%%  $$
 $$
 E_n(\widetilde\chi_h)_1\le\frac{1}{2nh}<1,\quad
  {h>\frac{1}{2n}}. \eqno (2)
  $$
}
%{In addition, we believe that (see \cite[section 5, figure on p.
%31]{bku}) $$
% E_n(\chi_h)_1 < \frac{1}{2nh},\quad h\in{(0,1]}\setminus M_n. \eqno (5'')
%$$ }

For $h\in{(0,1]}\setminus M_n$ we used
the precise description of the signum--functions from $C_n^\perp$.
Denote by  $G_n$ a class of functions $g(t)$ with the following
properties:
\begin{itemize}
%%%%%%%%    \item $g(t)=g(-t)$;  поскольку  $g\in C_n^\perp$
    \item $|g|=1${,}
    \item $g  \in C_n^\perp${,}
    \item
       function $g$ has  $ n+1 $ breakpoints on  $(0,{1/2})$:
    \ {$t_{0,n}<t_{1,n}<\ldots<t_{n,n}.$}
\end{itemize}

{The following Lemma has direct links to  results of P.Thchebyshev
(1859), A.Markov
 (1906), S.Bern\-stein(1912), Y.Ge\-ro\-nimus (1935),
 G.Szego (1964), F.Pe\-her\-s\-tor\-fer (1979)
 (see \cite{bku} for some details).}

{\bf Lemma A.} { \it The set of the zeros of the equation
$$
\cos 2\pi(n+1)t  - 2q \cos 2\pi nt + q^2 \cos 2\pi (n-1)t =0, \quad  q \in
(-1,1), \eqno (3)
$$
on  $(0, 1/2)$ is equal to the set of the breakpoints of some function from  $G_n$.
In the converse direction: for any function $g_0$ from  $G_n$ there is $q_0 \in (-1,1)$ such that
the set of zeros of {$(3)$} on  $(0, 1/2)$ is equal to the set of breakpoints of $g_0$.}
\vskip .2cm
\noindent
Lemma A {and $(1'')$ give} formula for the best approximations of characteristic function
for arbitrary {$h>0$} (\cite{bku}, Theorem 5). In particular,
the following statement {(see \cite[section~5,
p.~30]{bku})} is true.
\vskip .2cm
{\bf Theorem B.} For\; {$n\ge2$,}\; $h={t_{1,n}}\in(1/(2n),3/(2n))$
$$
E_n(\widetilde\chi_h)_1 = 1 - 2 t_{0,n}/t_{1,n}.
$$

\vspace*{-8mm}
\begin{center}
%\begin{figure}[htb]
%\includegraphics[height=8cm]{R_j2.eps}
\includegraphics[trim=2mm 2mm 2mm 2mm,clip,height=14cm,angle=270]{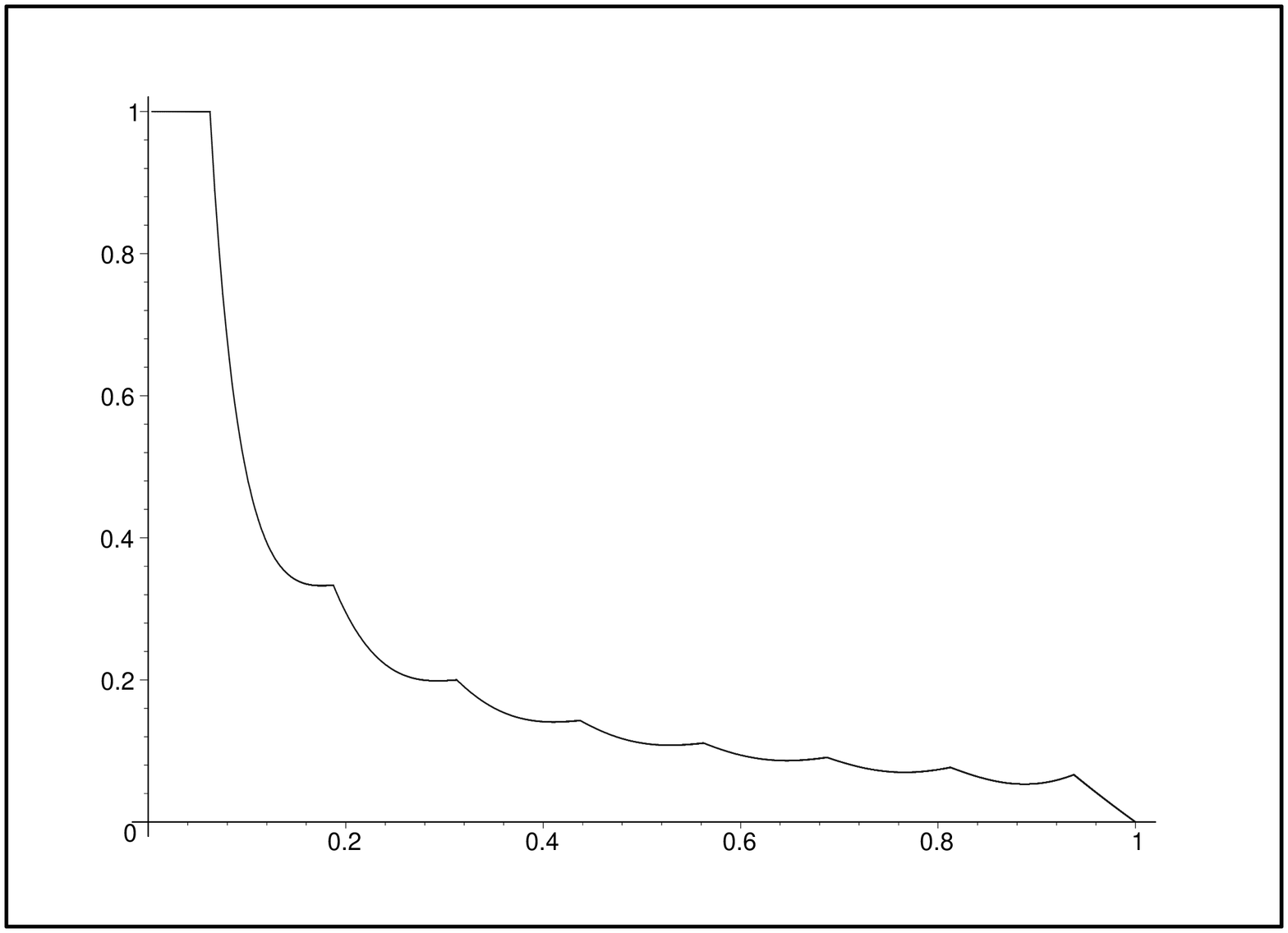}
%\caption[]

%\end{figure}

  {\footnotesize Plot of $\Psi_8(h):=E_8(\widetilde\chi_h)_1$
  for $h\in\left[0,1\right].$}\label{Ris-Psi-8}
\end{center}

\vspace*{-8mm}
\begin{center}
%\begin{figure}[htb]
%\includegraphics[height=8cm]{R_j2.eps}
\includegraphics[trim=2mm 2mm 2mm 2mm,clip,height=14cm,angle=270]{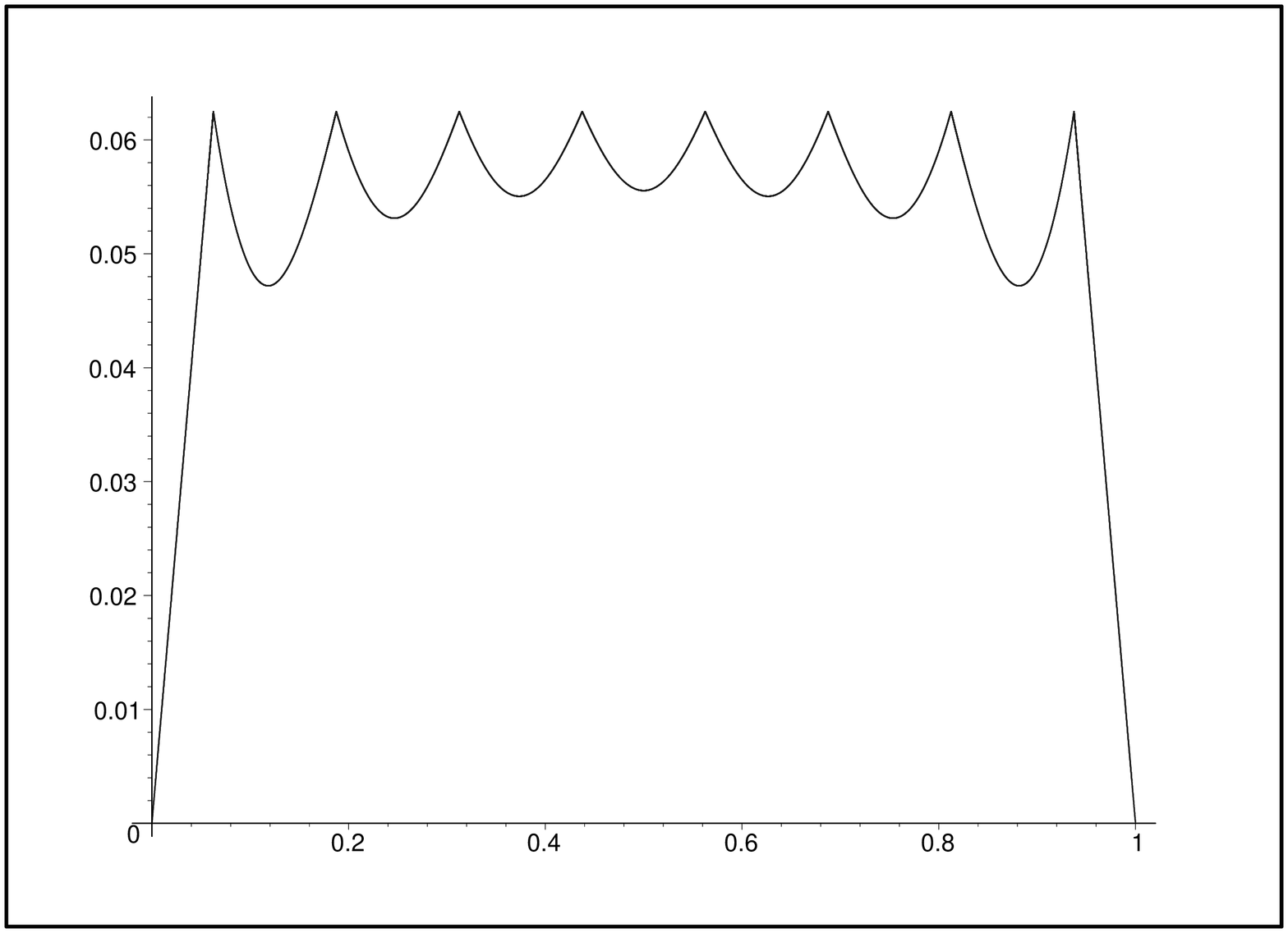}
%\caption[]

%\end{figure}

  {\footnotesize Plot of $\Phi_8(h):=h{\cdot}E_8(\widetilde\chi_h)_1$
  for $h\in\left[0,1\right].$}\label{Ris-Phi-8}
\end{center}

\vskip .5cm
\begin{center}
{\bf 4  Evaluation of  $E_n(\widetilde\chi_h)_1$ for $h={1/n}$}
\end{center}
\vskip .2cm

The case $h={1/n}$ is important for two reasons.
\vskip .2cm

{\bf 1.} This case is the  start point of our investigations on
approximation of concrete functions in $L$--metrics. The answer
to this question is the principal particular case of Theorem B.

{\bf 2.}  We have nice formula in this case {(see \cite[section~5,
p.~30]{bku}).}

\vskip .2cm

{\bf Theorem  C.}{ \it For $n \ge 2, \ h(n)=1/n$ we have $$ c(h(n), n) =
E_n(\widetilde\chi_{h(n)})_1 <  E_{n+1}(\widetilde\chi_{h(n+1)})_1 < \cdots  < \lim_{k
\to \infty} E_k(\widetilde\chi_{h(k)})_1 = 1- 2 v_0, $$ where $v_0$ is the
first positive zero of equation $$ \sec \pi v - \tan \pi v = v. $$
} Note that (see ($2$)) $$ 1-2 v_0 = 0.3817350529 \dots  < 1/2.
$$
\vskip .5cm
\begin{center}
{\bf 5  Favard and Jackson type theorems  }
\end{center}
\vskip .2cm

Put
$$
 W_2(f,h,x):= f(x) - (f*\chi_h )(x)=-\frac 1h
\int_0^{h/2} \left(  f(x-t)-2f(x)+f(x+t) \right) \, dt, \eqno (4)
$$ $$ W_2 (f,h):=\sup_x |f(x) - (f*\chi_h) (x)|. \eqno (5)
 $$
It is clear that
 $$
 2W_2(f,h)\le\omega_2(f,h/2):=\sup_{x,\,0<t<h}|f(x-t/2) - 2 f(x) + f(x+t/2)|. \eqno (6)
 $$

{\bf Theorem 1 (Favard type)}.{ \it Let $g\in T_{2n-1}^\perp $.
Then $$ \| g \| \le (1- c(h,n))^{-1} W_2 (g,h), \quad
{h>\frac{1}{2n}}. \eqno (7) $$ }

{\it Proof}. It is  { a direct consequence} of {$(4),(5),(1), {(1')},(2)$.} The identity
$$
g(x)=(g*\chi_h)(x)+ W_2(g,h,x)
$$
yields
$$
\|g\|\le\|g*\chi_h\|+ W_2(g,h)\le c(n,h)\|g\|+W_2(g,h).
$$
\qed

Denote a space of continuous functions on  $T$ by $C(T)$.

\medskip

%%%% This statement is equivalent to following equalities for Fourier
%%%% coefficients:
%%%%
%%%% $$
%%%% \widehat{\tau_k} = \widehat{g_k} \left( \widehat{(\delta - \chi_h
%%%% + \tau_h)_k} \right)^{-1}, \quad  k=-n+1,\dots, n-1.
%%%% $$

%\medskip

{\bf Theorem 2 (Jackson type).}{ \it Let $f \in C({T})$. Then for
$ {h>\frac{1}{2n}}$
$$
E_n(f) :=\inf_{\tau\in T_{2n-1}}
\|f-\tau\|\le(1-c(h,n))^{-1}W_2(f,h).\eqno (8)
$$
}
{\it Proof}. The inequality (8) is {a} modification of (7).
If $\tau_h$ is the
best  $L$--approximation of  $\widetilde\chi_h$, then for
suitable choose of $ \tau_f \in T_{2n-1} $ we have
$$
(f-\tau_f)(x) = (f-\tau_f)\odot(\widetilde\chi_h - \tau_h)(x) + W_2(f,h,x),
$$
and we can repeat the proof of Theorem 1. \qed

\vskip .5cm
\begin{center}
{\bf 6  On exact constants in Favard and Jackson theorems }
\end{center}
\vskip .2cm

The constants in Theorems 1 and 2 are not sharp.

Consider the principal case $h=(2n)^{-1}$. Theorems 1 and 2
take place in this case too, but for the proof we need more
complex ideas.
%considerations
%complex ideas
(see chapter 8).
\vskip .3cm
{\bf Conjecture}. {\it The following inequalities are true

$$
\| g \| \le 3 \,  W_2 (g,1/(2n)), \quad g \in T_{2n-1}^\perp,
\eqno (9)
$$
$$
E_n(f) \le 3 \, W_2 (f,1/(2n)), \quad f \in C(T).  \eqno
(10)
$$
}

We can not replace the constants 3 in the inequalities  $(9),
(10)$ with smaller constants. The  inequalities ($9$), ($10$) imply
the sharp Favard's inequalities as follows.

$$
\| g \| \le 3 W_2(g,(2n)^{-1}) \le 3 \cdot 2n \, {\|D^2 g\|} \int_0^{1/(4n)} t^2  \, dt
= 2^{-5} n^{-2} \|D^2 g \|, \quad g \in T_{2n-1}^\perp,
$$
{
$$
E_n(f) \le 3 W_2(f,(2n)^{-1}) \le 3 \cdot 2n \, {\|D^2 f\|} \int_0^{1/(4n)} t^2  \, dt
= 2^{-5} n^{-2} \|D^2 f \|, \quad f \in C^2(T).
$$
}

\vskip .5cm
\begin{center}
{\bf 7  On the classical Favard and Jackson inequalites }
\end{center}
\vskip .2cm
The inequalities  (6), (8) and
$$
W_2(f,h) \le \frac {\|D^2 f\|} {h} \int_0^{h/2} t^2  \, dt
=\frac {h^2}{24} \|D^2 f\|,
$$
give the classical Favard and Jackson inequalities for the second derivative and
the second modulus of smoothness (this means that it  is also  true for the first
modulus of continuity).

\vskip .5cm
\begin{center}
{\bf 8  On the extrapolation of Favard and Jackson inequalities}
\end{center}
\vskip .2cm

Despite the fact that approximation of the characteristic function
$\widetilde \chi_h$
is possible only from some value of support $h$ we can prove the inequality
$(7)$ and  $(8)$ for small values of $h$ too. Let us show how to do this.

Write the identity
%\Red{for $g\in L(T)$}
 $$
 g=g-g*\chi_h+\chi_h*(g-g*\chi_h)+ \chi_h^2*g=
 $$
 $$
 g - g*\chi_h + \widetilde\chi_h\odot(g-
 g*\chi_h) + \widetilde\chi_h^2 \odot g. \eqno (11)
 $$

Denote by $E_n(\widetilde\chi^2)_1$ the best
$L$--approximation of
$$ \widetilde \chi_h^2:=\widetilde \chi_h \odot \widetilde \chi_h $$ by
{$T_{2n-1}$}.

Then, {using slightly modified proof of Theorem 1 we get:}

$$
\|g\| \le (1 + 1) W_2 (g,h) + E_n(\widetilde\chi_h^2)\|g \|, \quad g \in T_{2n-1}^\perp
$$
\noindent
and

$$ \|g \| \le 2(1-E_n(\widetilde\chi_h^2))^{-1} W_2(g,h).  $$
 The
last inequality is valid for some  $h$ such that
 $$
E_n(\widetilde\chi_h^2) < 1.  $$
 It is known that  (see \cite{bka3} ) $$
E_n(\widetilde\chi_h^2)=1/2, \quad h=1/(2n). $$ So, we can take the constant
in the inequality  $(9)$  equal to $4$. We can use the father
extrapolation of  $(11)$.
{ Identity}

$$ g=g - g*\chi_h + \chi_h*(g- g*\chi_h) + \chi_h^2*(g-g*\chi_h)+
\dots   =
$$
$$
g - g*\chi_h + \widetilde \chi_h \odot (g- g*\chi_h) + \widetilde \chi_h^2 \odot(g-g*\chi_h)+
\dots
\eqno (12) $$

{gives}

$$
{ \| g \| \le \left(1+ \sum_{j=1}^\infty E_n (\widetilde\chi_h^j)_1 \right) \  W_2 (g,h).}
$$
\noindent
In the case $h=1/(2n)$ (see \cite{bka3})
$$
{ E_n(\widetilde\chi_h^j)_1 =  F_j,}
$$
where
$$
F_j= 2 \left(\frac{2}{\pi}\right)^{j+1}
\sum_{k=-\infty}^\infty (4k+1)^{-j-1}= \left(\frac 2\pi \right)^j \mathcal{K}_j.
$$
{In particular,}
$$
{ F_0=1,} \quad F_1=1, \quad F_2=1/2, \quad F_3=1/3, \quad F_4 =
5/24, \quad F_5= 2/15.
$$

\noindent
From
$$
\sum_{j=0}^\infty F_j = \sum_{j=0}^\infty \mathcal{K}_j (2/\pi)^j  =\sec (1) + \tan (1)
$$
\noindent
we deduce the estimate $(9)$ with the constant
 $
 \sec (1) + \tan (1) = 3.408223443\dots
 $
%\Red{(instead of 3)}.

The proof of Jackson's inequality with the same constant (Favard's constant) one can obtain in the
following way.
Let $\tau_h^j \in T_{2n-1} $ be the polynomial of the best $L$--approximation of $\widetilde\chi_h^j$.
Put
$$
\tau_{f,N} := \sum_{j=0}^{N-1} \tau_h^j \odot (f - f*\chi_h ) + \tau_h^{N}\odot f, \quad {\tau_h^0=0.} \eqno (13)
$$
Then by subtraction $(13)$ from
$$
f= \sum_{j=0}^{N-1} \widetilde\chi_h^j \odot(f-f*\chi_h)+ \widetilde\chi_h^N \odot f, \quad {\widetilde\chi_h^0=\widetilde\delta,} \eqno (14)
$$
we get Jackson's {type} theorem:
$$
\| f - \tau_{f,N} \| \le  \left(\sum_{j=0}^{N-1} E_n (\widetilde\chi_h^j)_1 \right) \  W_2
(f,h)+ E_n (\chi_h^N)_1 \, \|f \|,
$$
and
$$
E_n (f) \le \left(\sum_{j=0}^\infty E_n (\widetilde\chi_h^j)_1 \right) \  W_2 (f,h). \eqno
(15)
$$

\vskip .5cm
\begin{center}
{\bf 9  On Stechkin's theorem  }
\end{center}
\vskip .1cm

\begin{center}
{\bf 9.1  Introduction }
\end{center}
\vskip .1cm

Stechkin's theorem is the generalization of Jackson's theorem to differences of higher
orders:
$$
\Delta_t^r f(x):=\sum_{j=0}^k (-1)^{j} \binom rj f(x + j t).
$$
Classical Stechkin's inequality is formulated in notation of
$r$--th modulus of smoothness
$$
\omega_r (f,h):= \sup_{0<t\le h} \| \Delta_t^r f \|,
$$
and has the following form:
$$
E_n (f) \le K_{n,r} (h) \ \omega_r (f, h).
$$
The behavior of the sharp constant (as the function of $n,r,h$ )
$$
K_{n,r}(h):=\sup_{f \in C} \frac{E_n(f)}{\omega_r (f,h)}
$$
is not clear in details.
Put
$$
\gamma_r^*:=  \begin{cases}
{\binom {2k}k}^{-1} {,} \quad & r=2k {,} \\
{\binom {2k-1}{k-1}}^{-1}{,} \quad & r=2k-1
\end{cases}\qquad
{\left(\gamma_r^*\asymp \frac{r^{1/2}}{2^r}\right)}.
$$
\noindent
It was recently proved  (see \cite{fks}) that
$$
K_{n,r}(\alpha/(2n))  \le C_\alpha \gamma_r^* , \quad
\alpha >1, \eqno (16)
$$
$$
K_{n,r}(1/(2n)) \le C \sqrt{r} \ln (r+1) \gamma_r^*. \eqno (17)
$$
In particular, for $\alpha = 2$

$$
K_{n,r}(1/n)  \le  5 \gamma_r^*.
$$
\noindent
Inequality  $(16)$  is not true for  $ \alpha < 1$ \cite{fks}. So, we have the
intrinsic open question: is the inequality

$$
K_{n,r}(1/(2n) )  \le  C  \gamma_r^*.
$$
\noindent
true?

We will show that the method of chapter 8 {allows} us to prove that

$$
K_{n,r}(1/n)  \le  \sqrt{2} \gamma_r^*, \quad K_{n,r}(1/(2n)) \le C \sqrt{r} \, \gamma_r^*.  \eqno (18)
$$
\noindent
It is known (\cite{fks}) that for  $h \le 1/(2k)$
$$
K_{n,r}(h) \ge c' \gamma_r^*,
$$
where
$$
c':= \begin{cases} \frac r{r+1}, \quad & r=2k-1; \\
1, \quad & r=2k.
\end{cases}
$$
Therefore, in the classical case  $\delta = 1/n$ we have the narrow interval for
the value $K_{n,r}(\delta)$.

\vskip .5cm

\begin{center}
{\bf 9.2 Smoothness and general results   }
\end{center}
\vskip .2cm

We assume that the smoothness order is {an} even number. In other
words, we suppose that $r=2k$.
It is convenient  to consider the symmetric differences:
$$
\widehat \Delta_t^{2k} f(x) = \sum_{j=-k}^{k} (-1)^j \binom {2k}{k+j} f
(x+jt).
$$
Introduce a class of even, integrable {functions}
$\Phi$. This is a class of the convolution kernels.
We write  $\phi \in \Phi$ if $\phi$ {is integrable on $R$ function with compact
support and}
$$
\phi (x) = \phi (-x), \quad  \int_{{{R}}} \phi (t) \, dt =1.
$$
We will use notation
$$
\phi_j (x):= \frac 1j \phi \left( \frac xj \right), \quad \phi_0 (t)
:=\delta (x).
$$
Define the function, measuring the $2k$--th  $\phi$--th smoothness of $f$ at the point $x$.
$$
W_{2k} (f,\phi, x):= {\binom {2k}k}^{-1} \int_{{{R}}}
\widehat \Delta_t^{2k} f(x) \phi (t) \, dt.
$$
The function  $W_{2k} (f,\phi, x)$ {can be written} as the convolution of  $f$
with the function
$$
W(x):=W_{2k} (\phi,x):= {\binom {2k}k}^{-1} \sum_{j=-k}^k (-1)^j
\binom{2k}{k+j} \phi_j (x) =
$$
$$
\delta (x) - 2 \sum_{j=1}^k
(-1)^{j+1} a_j \phi_j(x), \quad
a_j:=\frac{\binom{2k}{k+j}}{\binom{2k}{k}}.
$$
Put
$$
U(x):=U_{2k} (f,\phi, x):= 2 \sum_{j=1}^k (-1)^{j+1} a_j \phi_j(x).
$$
and denote by  $U^j$  the convolution power of $U$:
$$
U^0(x):= \delta (x), \quad U^j(x):=(U*U^{j-1})(x).
$$
\noindent
The identity  ( see (12))

$$ g = \sum_{j=0}^\infty U^j*(g-U*g) = \sum_{j=0}^\infty
\widetilde U^j \odot (g-U*g) $$ and equalities $$ W_{2k}(f,\phi,x)
=  (W*f)(x) = (f-U*f)(x) $$ give the following result \vskip .2cm
{\bf Theorem  3 (Favard type)}.{ \it Let $g \in T_{2n-1}^\perp $.
Then $$ \| g \| \le  \left( \sum_{j=0}^\infty E_n (\widetilde
U^j)_1 \right) \| W_{2k} (g,\phi, \cdot  ) \|. $$ } \noindent The
passage from Theorem 3 to Stechkin's inequality {is described} in
the last lines of
 chapter 8 (see $(13)-(15)$).

\vskip .3cm

{\bf Theorem 4 (Stechkin type)}.{ \it Let $f \in C(T) $. Then
$$
E_n(f) \le  \left( \sum_{j=0}^\infty E_n (\widetilde U^j)_1 \right) \| W_{2k} (f,\phi, \cdot  ) \|.
$$
}

\begin{center}
{\bf 9.3  Concrete results  }
\end{center}
\vskip .2cm
\noindent
By choosing
$$
\phi (x) = \chi_h^2 (x)
$$
one can obtain the estimate \cite{fks}
$$
\sum_{j=0}^\infty E_n (\widetilde U^j)_1
\le \left(  \cos \frac{\pi}{2} \rho  \right)^{-1},
\quad \rho = \frac{ \mu_{2k}}{2 n h } < 1,  \quad \mu_{2k} \approx  ( 1 -
(2k)^{-1/2})^{1/2}.
$$
\noindent
Since
$$
 \| W_{2k} (f,\chi_h^2 , \cdot  ) \| \le \gamma_{2k}^* \omega_{2k}
 (f,h),
$$
 we have
\vskip .3cm
{\bf Theorem 5}.{ \it Let $f \in C(T), \  \alpha > 1 $. Then
$$
E_n(f) \le  (\cos (\pi/(2\alpha) )^{-1} \ \gamma_{2k}^*  \  \omega_{2k} \left(f, \frac{\alpha }{2n}  \right).
$$
}
\noindent
In particular, for  $\alpha = 2$
$$
(\cos (\pi/(2\alpha) )^{-1} = \sqrt{2}.
$$
\noindent
From
$$
\gamma_{2k}^* \omega_{2k} (f,\delta) \le \gamma_{2k-1}^*
\omega_{2k-1}(f,\delta),
$$
the first inequality (18) follows.

Consider the case $\alpha =1$. {Choose  $\phi = \chi_h^2$. We
can }
obtain good estimates in this case, however only
in terms of the characteristic
$W_{2k}(f,\chi_h^2,x)$.
\noindent
Put
$$
K_{n,2k}^*(\phi):=\sup_{f \in C} \frac{E_n(f)}{\| W_{2k} (f,\phi, \cdot
)\|}.
$$
\vskip .3cm
\noindent
We can present now the corollary of  Theorem 7.1 from
\cite{fks}.
\vskip .3cm

{\bf Theorem 6}. { For $r=2k$ we have
$$
\frac{\gamma_{r}^*}{1-\mu_{r}^2}  \le  K_{n,r}^*(\chi_{1/(2n)}^2) \le \frac 4 \pi
 \  \frac{\gamma_{r}^*}{1-\mu_{r}^2},
$$
{where}
$$
\frac{\gamma_{r}^*}{1-\mu_{r}^2} \asymp \sqrt{r} \gamma_{r}^*
\asymp \frac{r}{2^r}.
$$
}
So we can omit the factor $\ln (r+1)$ in  $(17)$.

\vskip .2cm
\vskip .5cm

\begin{center}
{\bf 10 Comments   }
\end{center}
\vskip .2cm

\vskip .2cm
\noindent
{\bf 1. Equivalence of convolutions}

Denote by $\widetilde{g}$ the 1-periodization of  $g {\in L(R),
\quad \mbox{supp} \, g < \infty }$ :
 $$
 \widetilde{g}(t):=\sum_{k\in Z}g(t+k).
 $$
Let $f\in L({T})$ be an arbitrary 1-periodic function. It is
{well--known } that
 $$
 (f\odot{\widetilde g})(x):=\int_{T}f(x-t)\widetilde{g}(t)\,dt=
 \int_{-\infty}^{+\infty}f(x-t){g}(t)\,dt =: (f * g) (x).
 $$
Indeed,
 $$
 I(x):=\int_{T}f(x-t)\widetilde{g}(t)\,dt=\int_{0}^{1}f(x-t)\sum_{k\in Z}g(t+k)\,dt
 =\sum_{k\in Z}\int_{0}^{1}f(x-t)g(t+k)\,dt.
 $$
The change of variable $t+k=u,$ $t=u-k$ gives
 $$
 I(x)=\sum_{k\in Z}\int_{0}^{1}f(x-t)g(t+k)\,dt=
 \sum_{k\in Z}\int_{k}^{k+1}f(x-u+k)g(u)\,du=
$$
$$
 \sum_{k\in Z}\int_{k}^{k+1}f(x-u)g(u)\,du=
 \int_{-\infty}^{+\infty}f(x-t){g}(t)\,dt.\qed
 $$

\vskip .5cm
{\bf 2. How to choose $\tau_f$ in Theorem 2}

{\bf Lemma 1.} {\it Let $g\in L(T)$ and
 $$
 |\widehat{g}(k)|<1\quad\mbox{for all}\quad k\in Z,\quad\mbox{where}\quad
 \widehat{g}(k):=\int_{T}g(t)e^{-2\pi{i}k t}\,dt. \eqno (1c)
 $$
Then for arbitrary $n\in{N},$ $\varphi\in T_{2n-1}$ there exist $\tau\in T_{2n-1}$ such
that}
 $$
 \varphi=\tau-\tau \odot g. \eqno(2c)
 $$
{\it Proof}. The equation $(2c) $ is equivalent to
 $$
 \varphi=\tau \odot (D_n-g),\quad\mbox{where}\quad
 D_n({x})=\sum_{k=-(n-1)}^{n-1}e^{2\pi{i}k{x}},
 $$
and
 $$
 \widehat{\varphi}(k)=\widehat{\tau}(k)(1-\widehat{g}(k)),\quad |k|\le n-1.
 $$
Thus
 $$
 \tau({x})=\sum_{k=-(n-1)}^{n-1}\frac{\widehat{\varphi}(k)}
 {1-\widehat{g}(k)}e^{2\pi{i}k {x}}.
 $$

{\bf Remark 1.}  Condition $\|g\|_L<1$ imply $(1c)$.

\medskip
{\bf Remark 2.} We proved Theorem 2 in the following form:
$$
\|f-\tau_f\|\le\frac{W_2(f,h)}{1-c(h,n)},\quad f\in C({T}),\quad
n\in N,\quad \frac{1}{2n}<h\le1, \eqno (10') $$ where
 $$
 \tau_f(x)=\sum_{k=-(n-1)}^{n-1}
 \frac{\widehat{f}(k)\widehat{\tau}_h(k)}
 {1-\widehat{\chi}_h(k)+\widehat{\tau}_h(k)}e^{2\pi{i}kx}.
 $$

\vskip .5cm
{\bf 3. About the equality} $ E_n(\widetilde\chi_{1/(2n)}^j)_1 =  F_j.$
\vskip .2cm
In \cite{bka3} we proved this equality with some restrictions on $h$.
Equivalence of convolutions $*$ and $\odot$ for periodic functions
gives the proof without restrictions.
Namely, we do not need to modify anything in \cite{bka3}.

\vskip .5cm
{\bf 4. About the  inequality} $ \sum_{j=0}^\infty E_n (\widetilde
U^j)_1
\le \left(  \cos \frac{\pi}{2} \rho  \right)^{-1}, \quad \rho <1
$.
\vskip .2cm
In \cite{fks} inequality had been proved  in another form
$$
\sum_{j=0}^\infty \|  U^j \|_{T^\perp_{2n-1}}
\le \left(  \cos \frac{\pi}{2} \rho  \right)^{-1}, \quad \rho <1.
$$

The equality $$ \|  U^j \|_{T^\perp_{2n-1}} = E_n (\widetilde
U^j)_1 $$ allows us to simplify the approach to Stechkin's theorem
and give new estimates of constants. Namely, in this paper we
proved that Stechkin's constants (in Theorem 4) are equal to
Favard's constants (in Theorem 3).
\vskip .3cm
{\bf Acknowledgements.} Our thanks to Alexei Solyanik for his
comments on a draft of this paper.

\vskip .2cm

\small
\begin{tabular}[t]{l}
A.G. Babenko \\ Institute of Mathematics and Mechanics \\ Ural
Branch of the Russian Academy of Sciences \\ 16, S.Kovalevskoi
Str. \\ Ekaterinburg, 620219 \\ Russia \\
\textrm{babenko@imm.uran.ru}
\end{tabular}
\hfill
\begin{tabular}[t]{l}
Yuriy Kryakin \\
Institute of Mathematics \\
University of Wroclaw \\
Plac Grunwaldzki 2/4 \\
50-384 Wroclaw \\
Poland \\
\textrm{kryakin@math.uni.wroc.pl}
\end{tabular}

\begin{thebibliography}{1}


\bibitem{bk1} Babenko A.G., Kryakin Yu.V. \emph{On approximation of step functions
by trigonometric polynomials in the integral metric }// Izvestia
of the Tula State University. Ser. Mathematics. Mechanics.
Informatics. Tula: TSU, 2006, Vol. 12, N 1, 27--56.


\bibitem{bku} Babenko A.G., Kryakin Yu.V. \emph{Integral
approximation of characteristic function of interval by
trigonometric polynomials}// Trudy Instituta Matematiki i
Mekhaniki, 2008, Vol. 14, N 3, 19--37.


\bibitem{bka3}
Babenko A.G., Kryakin Yu.V.  \emph{On L--Approximation of
B--splines by trigonometric polynomials}
//ArXiv:0811.0686v1 [math. CA], (2008), 1--6.

\bibitem{fks} Foucart S., Kryakin Yu., Shadrin A. \emph{On the exact constant
in Jackson-Stechkin inequality for the uniform metric} //
ArXiv:math CA/0612283, (2006), 1--20 (to appear in Constructive
Approximation).



\bibitem{M} Markov A.A.,  \emph{Selected works} // Moscow, 1948,
146--230 (Russian).
\end{thebibliography}
\end{document}